\documentclass[12pt,draft]{amsart}
\usepackage{amsmath,amsthm,latexsym,amscd,amsbsy,amssymb,pb-diagram}
\chardef\bslash=`\\ 

\makeatletter
\def\verbatim{\interlinepenalty\@M \@verbatim
  \leftskip\@totalleftmargin\advance\leftskip2pc
  \frenchspacing\@vobeyspaces \@xverbatim}
\makeatother
\newcounter{rmnum}

\newtheorem{thm}{Theorem}[section]
\newtheorem{cor}[thm]{Corollary}

\newtheorem{pro}[thm]{Proposition}

\theoremstyle{definition}

\theoremstyle{remark}

\numberwithin{equation}{section}

\begin{document}


\title[Complemented subspaces of products of Banach spaces]
{Complemented subspaces of products of Banach spaces}
\author{Alex Chigogidze}
\address{Department of Mathematics and Statistics,
University of Saskatche\-wan,
McLean Hall, 106 Wiggins Road, Saskatoon, SK, S7N 5E6,
Canada}
\email{chigogid@math.usask.ca}
\thanks{Author was partially supported by NSERC research grant.}



\begin{abstract}{We show that complemented subspaces of
uncountable products of Banach spaces are products of complemented
subspaces of countable subproducts.} 
\end{abstract}

\maketitle
\markboth{A.~Chigogidze}{Complemented subspaces of products
of Banach spaces}

\section{Introduction}\label{S:intro}

The following old unsolved problem (L.~Nachbin \cite{nac61}) of
describing injective
locally convex spaces is one of the general problems of the structure
theory of locally  convex spaces.
\medskip

{\bf Problem 1}. Is every injective
locally convex space isomorphic to a product of Banach injective spaces?
\medskip

In investigations related to this problem (see, for instance,
\cite{domort89},
\cite{dom901}, \cite{dom902}, \cite{metmos891}) the following problem
(\cite[p.~71]{dom902}, \cite[p.~147]{metmos90}) arose.
\medskip

{\bf Problem 2}.
Is every complemented subspace of a product of a (countable) family of
Banach spaces isomorphic to a product of Banach spaces? 
\medskip

G.~Metafune and V.~B.~Moscatelli \cite[p.~251]{metmos892} conjectured that
this is false in general. Later this conjecture has been confirmed
by M.~Ostrovskii \cite{ost96} who showed that not all complemented
subspaces of countable
products of Banach spaces are isomorphic to products of Banach spaces.

Our main result shows that for uncountable products situation is somewhat
different.

{\bf Theorem.}
{\em A complemented subspace of an uncountable pro\-duct of Banach spaces
is a product of complemented subspaces of countable subproducts.}

The following immediate corollary of this result provides a partial solution
to Problem 1.

{\bf Corollary.}
{\em Every injective locally convex space is isomorphic to a product
of injective Fr\'{e}chet spaces.}

Author is grateful to P.~Doma\'{n}skii for his comments on the Internet
version of this note.


\section{Results}

The following statement expresses a key fact used in the proof of Theorem \ref{T:complementedpro}. 

\begin{pro}\label{P:spectral}
Let $\displaystyle r \colon \prod\{ B_{t} \colon t \in T \}
\to \prod\{ B_{t} \colon t \in T\}$
be a continuous linear map of an uncountable product
of Banach spaces into itself.
Let also $A$ be a countable subset of $T$. Then there
exist a countable subset $S$
of $T$ and a continuous linear map $\displaystyle r_{S} \colon
\prod\{ B_{t} \colon t \in S\} \to \prod\{ B_{t} \colon t \in S\}$
such that $A \subseteq S$ and $\pi_{S}\circ r = r_{S}\circ \pi_{S}$,
where $\displaystyle \pi_{S} \colon \prod\{ B_{t} \colon t \in T\}
\to \prod\{ B_{t} \colon t \in S\}$ denotes the projection
onto the corresponding subproduct.
\end{pro}
\begin{proof}
Let $\exp_{\omega}T$ denote the set of all countable subsets
of the indexing set $T$.
Consider the following relation

\begin{multline*}
 {\mathcal L} = \{ (S,R) \in \left( \exp_{\omega}T\right)^{2}
\colon  S \subseteq R\;\;\text{and there exists a continuous linear map}\\
 r_{S}^{R} \colon \prod\{ B_{t} \colon t \in R\} \to
\prod\{ B_{t} \colon t \in S\}\; \;\text{such that}\;\;
\pi_{S}\circ r = r_{S}^{R} \circ \pi_{R} \} ,
\end{multline*}

\noindent where 
\[ \pi_{S} \colon \prod\{ B_{t} \colon t \in T\} \to
 \prod\{ B_{t} \colon t \in S\}\]

\noindent and

\[ \pi_{S}^{R} \colon \prod\{ B_{t} \colon t \in R\} \to
\prod\{ B_{t} \colon t \in S\} \]

\noindent denote canonical projections onto the corresponding
subproducts.

We need to verify the following three properties of the above
defined relation.
\medskip

{\em Existence}. If $S \in \exp_{\omega}T$, then there
exists $R \in \exp_{\omega}T$
such that $(S,R) \in {\mathcal L}$.
\medskip

{\em Proof}. Let $S = \{ t_{n} \colon n \in \omega\}$.
For each $n \in \omega$ consider the composition
$\displaystyle \pi_{t_{n}} \circ r \colon
\prod\{ B_{t} \colon t \in T\} \to B_{t_{n}}$. Since
$B_{t_{n}}$ is a Banach space, it follows that every
continuous linear map into $B_{t_{n}}$, defined on an infinite product
of Banach space, can be factored through a finite subproduct
(this is a well known fact; see, for instance,
\cite[Proposition 0.1.9]{hel93}). Consequently there exist a
finite subset $R_{n}$ and a continuous linear map
$\displaystyle r^{R_{n}}_{t_{n}} \colon \prod\{ B_{t} \colon
t \in R_{n}\} \to B_{t_{n}}$ such that
such that $\pi_{t_{n}}\circ r = r^{R_{n}}_{t_{n}}\circ \pi_{R_{n}}$
for each $n \in \omega$. Without loss of generality we may
assume that $t_{n} \in R_{n}$ for each $n \in \omega$
(otherwise consider the set $R_{n} \cup \{ t_{n}\}$).
Let $R = \cup\{ R_{n} \colon n \in \omega\}$ and
$r^{R}_{t_{n}} = r^{R_{n}}_{t_{n}}\circ
\pi_{R_{n}}^{R}$, $n \in \omega$. Clearly
\[ \pi_{t_{n}}\circ r = r^{R_{n}}_{t_{n}}\circ
\pi_{R_{n}} = r^{R_{n}}_{t_{n}}\circ \pi_{R_{n}}^{R}
\circ \pi_{R} = r^{R}_{t{n}}\circ \pi_{R} ,\; n \in \omega .\]

Next consider the diagonal product
\[ r_{S}^{R} = \triangle\{ r^{R}_{t_{n}} \colon n \in \omega\}
\colon \prod\{ B_{t} \colon t \in R\} \to \prod\{ B_{t_{n}}
\colon n \in \omega\} = \prod\{ B_{t} \colon t \in S\} \]

\noindent and note that $r_{S}^{R}$ is a continuous linear map
which satisfies the equality
$\pi_{S}\circ r = r_{S}^{R}\circ \pi_{R}$. This shows
that $(S,R) \in {\mathcal L}$.
\medskip

{\em Majorantness}. If $(S,R) \in {\mathcal L}$,
$P \in \exp_{\omega}T$ and $R \subseteq P$,
then $(S,P) \in {\mathcal L}$.
\medskip

{\em Proof}. This is trivial. Indeed let
$\displaystyle r_{S}^{R} \colon \prod\{ B_{t} \colon t \in R\} \to
\prod\{ B_{t} \colon t \in S\}$ be a continuous linear map
such that $\pi_{S} \circ r = r_{S}^{R} \circ \pi_{R}$. Consider the map
$\displaystyle r_{S}^{P} \colon \prod\{ B_{t} \colon t \in P\} \to
\prod\{ B_{t} \colon t \in S\}$ defined as the composition
$r_{S}^{P} = r_{S}^{R}\circ \pi_{R}^{P}$. Since
$\pi_{S}\circ r = r_{S}^{R}\circ \pi_{R} =
r_{S}^{R}\circ \pi_{R}^{P}\circ \pi_{P} = r_{S}^{P}
\circ \pi_{P}$ it follows that $(S,P) \in {\mathcal L}$.

\medskip

{\em $\omega$-closeness}. Suppose that $(S_{i},R) \in
{\mathcal L}$ and $S_{i} \subseteq S_{i+1}$ for each
$i \in \omega$. Then $(\cup\{ S_{i} \colon
i \in \omega\}, R) \in {\mathcal L}$.
\medskip

{\em Proof}. Consider the following projective sequence

\[ \prod\{ B_{t} \colon t \in S_{0}\}
\xleftarrow{\pi_{S_{0}}^{S_{1}}}\cdots\xleftarrow{ }
\prod\{ B_{t} \colon t \in S_{i}\}
\xleftarrow{\pi_{S_{i}}^{S_{i+1}}}\prod\{ B_{t}
\colon t \in S_{i+1}\} \xleftarrow{ }\cdots\]
\medskip

\noindent limit of which is isomorphic to the product
$\displaystyle \prod\{ B_{t} \colon t \in S\}$,
where $S = \cup\{ S_{i} \colon i \in \omega\}$.

Since $(S_{i},R) \in {\mathcal L}$, there exists a continuous
linear map $\displaystyle r_{S_{i}}^{R} \colon
\prod\{ B_{t} \colon t \in R\} \to \prod\{ B_{t} \colon
t \in S_{i}\}$ such that $\pi_{S_{i}}\circ
r = r_{S_{i}}^{R}\circ \pi_{R}$, $i \in \omega$. Note that
$\pi_{S_{i}}^{S_{i+1}}\circ r_{S_{i+1}}^{R} = r_{S_{i}}^{R}$
for each $i \in \omega$.
Indeed let $\displaystyle x \in \prod\{ B_{t} \colon t \in R\}$
and consider any point $\displaystyle y \in \prod\{ B_{t} \colon t \in T\}$
such that $x = \pi_{R}(y)$. Since $(S_{i},R),
(S_{i+1},R) \in {\mathcal L}$ we have
\begin{multline*}
 \pi_{S_{i}}^{S_{i+1}}\left( r_{S_{i+1}}^{R}(x)\right) =
\pi_{S_{i}}^{S_{i+1}}\left( r_{S_{i+1}}^{R}
\left(\pi_{R}(y)\right)\right) = \pi_{S_{i}}^{S_{i+1}}
\left( \pi_{S_{i+1}}\left( r(y)\right)\right) =
\pi_{S_{i}}\left( r(y)\right) =\\
 r_{S_{i}}^{R}\left(\pi_{R}(y)\right) = \pi_{S_{i}}^{R}(x) .
\end{multline*}

\noindent In this situation the collection $\displaystyle
\left\{ \pi_{S_{i}}^{R} \colon \prod\{ B_{t} \colon t \in R\}
\to \prod\{ B_{t} \colon t \in S_{i} \} \colon i \in \omega\right\}$

\noindent uniquely defines a continuous linear map
$\displaystyle r_{S}^{R} \colon \prod\{ B_{t} \colon t \in R\}
\to \prod\{ B_{t} \colon t \in S\}$ such that
$\pi_{S_{i}}^{S}\circ r_{S}^{R} = r_{S_{i}}^{R}$ for each $i \in \omega$
($r_{S}^{R}$ is simply the diagonal product of $r_{S_{i}}^{R}$'s).
It only remains to note that $\pi_{S}\circ r = r_{S}^{R}\circ \pi_{R}$ which
completes the proof of the fact that $(S,R) \in {\mathcal L}$.

According to \cite[Proposition 1.1.29]{chibook96} the set of
${\mathcal L}$-reflexive elements of $\exp_{\omega}T$ is
cofinal in $\exp_{\omega}T$. An element $S \in \exp_{\omega}T$
is ${\mathcal L}$-reflexive if $(S,S) \in {\mathcal L}$.
In our situation
this means that the given countable subset $A$ of $T$ is
contained in a larger countable subset $S$ for which there
exists a continuous linear map
$\displaystyle r_{S} = r_{S}^{S} \colon \prod\{ B_{t}
\colon t \in S\} \to \prod\{ B_{t} \colon t \in S\}$ satisfying the equality
$\pi_{S}\circ r = r_{S} \circ \pi_{S}$. Proof is completed.
\end{proof}

\begin{thm}\label{T:complementedpro}
A complemented subspace of a product of uncountable family
of Banach spaces is isomorphic to a product
of Fr\'{e}chet spaces.
More formally, if $X$ is a complemented subspace of the product
$\displaystyle \prod\{ B_{t} \colon t \in T\}$ of
Banach spaces $B_{t}$, $t \in T$, then $X$ is
isomorphic to the product $\displaystyle \prod\{ F_{j} \colon j \in J\}$,
where $F_{j}$ is a complemented subspace of the product
$\displaystyle \prod\{ B_{t} \colon t \in T_{j}\}$ with $|T_{j}| = \omega$
for each $j \in J$.
\end{thm}
\begin{proof}
Let us first of all set up a notation. For a subset $S \subseteq T$,
where $T$ is an indexing set with $|T| = \tau >\omega$, let
\[ B_{S} = \prod\{ B_{t} \colon t \in S\} \;\; \text{and}\;\;
B = \prod\{ B_{t} \colon t \in T\} . \]

Let also for $S \subseteq R \subseteq T$ 

\[ \pi_{S} \colon B = \prod\{ B_{t} \colon t \in T\} \to
B_{S} = \prod\{ B_{t} \colon t \in S\}\]

\noindent and

\[ \pi_{S}^{R} \colon B_{R} = \prod\{ B_{t} \colon t \in R\} \to
B_{S} = \prod\{ B_{t} \colon t \in S\} \]

\noindent denote canonical projections onto the corresponding subproducts.

Let $X$ be a complemented subspace of the product
$\displaystyle B = \prod\{ B_{t} \colon t \in T\}$.
Choose a continuous
homomorphism $r \colon B \to X$ such that $r(x) = x$ for each
$x \in X$. Let us agree that a subset $S \subseteq T$ is called
$r$-admissible if
$\pi_{S}\left( r(z)\right) = \pi_{S}(z)$ for each point
$z \in \pi_{S}^{-1}\left(\pi_{S}(X)\right)$. 
\medskip

{\bf Claim 1}. {\em The union of an arbitrary family of
$r$-admissible sets is $r$-admissible.}
\medskip

Let $\{ S_{j} \colon j \in J\}$ be a collection of
$r$-admissible
sets and $S = \bigcup\{ S_{j} \colon j \in J\}$.
Let $z \in \pi_{S}^{-1}\left( \pi_{S}(X)\right)$.
Clearly $z \in \pi_{S_{j}}^{-1}\left(\pi_{S_{j}}(X)\right)$
for each $j \in J$ and consequently
$\pi_{S_{j}}\left( r(z)\right) = \pi_{S_{j}}(z)$
for each $j \in J$.
Assuming that there is a point
$z_{0} \in \pi_{S}^{-1}\left(\pi_{S}(X)\right)$ such that
$\pi_{S}(r(z_{0})) \neq \pi_{S}(z_{0})$ we conclude that
there exists an index $s \in S$ such that
$\pi_{\{ s\}}^{S}\left(\pi_{S}(r(z_{0}))\right) \neq \pi_{\{ s\}}^{S}\left(\pi_{S}(z_{0})\right)$. Since
$S = \bigcup\{ S_{j} \colon j \in J\}$
it follows that there exists an index $j \in J$ such
that $s \in S_{j}$. Then we have
$\pi_{ S_{j}}^{S}\left(\pi_{S}(r(z_{0}))\right) \neq \pi_{ S_{j}}^{S}\left(\pi_{S}(z_{0})\right)$. But this is impossible

\[ \pi_{S_{j}}^{S}\left(
\pi_{S}\left( r(z_{0})\right)\right) = \pi_{S_{j}}
\left( r(z)\right) = \pi_{S_{j}}(z) =
\pi_{S_{j}}^{S}\left(\pi_{S}(z_{0})\right) .\]

\noindent This contradiction proves the claim.
\medskip

{\bf Claim 2.} {\em If $S \subseteq T$ is $r$-admissible, then
$\pi_{S}(X)$ is a complemented subspace of
$\displaystyle B_{S} = \prod\{ B_{t} \colon t \in S\}$.}
\medskip

Indeed, let $i_{S} \colon B_{S} \to B$ be the
canonical section of $\pi_{S}$
(this means that $i_{S} = \operatorname{id}_{B_{S}}\triangle
\mathbf{0} \colon B_{S} \to B_{S} \times B_{T-S} = B$).
Consider a continuous linear map 
$r_{S} = \pi_{S}\circ r\circ i_{S} \colon B_{S} \to \pi_{S}(X)$.
Obviously, $i_{S}(y) \in \pi_{S}^{-1}\left(\pi_{S}(X)\right)$
for any point $y \in \pi_{S}(X)$. Since $S$ is $r$-admissible
the latter implies that
\[ y = \pi_{S}\left(i_{S}(y)\right) =
\pi_{S}\left( r\left(i_{S}(y)\right)\right) = r_{S}(y) .\]
This shows that $\pi_{S}(X)$ is a complemented subspace of $B_{S}$.
\medskip

{\bf Claim 3}. {\em Let $S$ and $R$ be $r$-admissible subsets
of $T$ and $S \subseteq R \subseteq T$. Then 
there exists a topological isomorphism
$h_{S}^{R} \colon X_{R} \to X_{S} \times
\ker\left( \pi_{S}^{R}\right)$ which makes the diagram
\[
\begin{diagram}
\node{\pi_{R}(X)} \arrow{s,l}{\pi_{S}^{R}|\pi_{R}(X)}\arrow{e,t}{h_{S}^{R}}
\node{X_{S} \times \ker\left( \pi_{S}^{R}|\pi_{R}(X)\right)} \arrow{sw,b}{pi_{1}}\\
\node{X_{S}} 
\end{diagram}
\]

\noindent commutative.}
\medskip

Obviously $\pi_{R}(X) \subseteq \pi_{S}(X) \times B_{R-S}
\subseteq B_{R} = B_{S} \times B_{R-S}$. Consider the map
$i_{R} = \operatorname{id}_{B_{R}}\triangle
\mathbf{0} \colon B_{R} \to B_{R} \times B_{T-R} = B$. Also
let $r_{R} = \pi_{R}\circ r\circ i_{R} \colon B_{R} \to \pi_{R}(X)$.

Observe that
$\pi_{S}^{R}\circ r_{R}|\left( \pi_{S}(X) \times B_{R-S}\right) =
\pi_{S}^{R}|\left( \pi_{S}(X) \times B_{R-S}\right)$. Indeed, if
$x \in  \pi_{S}(X) \times B_{R-S}$, then $i_{R}(x)
\in \pi_{S}^{-1}\left(\pi_{S}(X)\right)$. Since $S$ is $r$-admissible,
we have $\pi_{S}\left(r\left(i_{R}(x)\right)\right) =
\pi_{S}\left(i_{R}(x)\right)$. Consequently,
\begin{multline*}
\pi_{S}^{R}\left( r_{R}(x)\right) =
\pi_{S}^{R}\left(\pi_{R}\left( r\left(i_{R}(x)\right)\right)\right) = \pi_{S}\left(r\left(i_{R}(x)\right)\right) =
\pi_{S}\left(i_{R}(x)\right) =\\
\pi_{S}^{R}\left(\pi_{R}\left( i_{R}(x)\right)\right) =
\pi_{S}^{R}(x) .
\end{multline*}

Next observe that $r_{R}(x) = x$ for any point $x \in \pi_{R}(X)$.
Indeed, since $R$ is $r$-admissible and since
$i_{R}(x) \in \pi_{R}^{-1}\left(\pi_{R}(X)\right)$ we have
\[ r_{R}(x) = \pi_{R}\left(r\left(i_{R}(x)\right)\right) =
\pi_{R}\left(i_{R}(x)\right) = x .
\]

In this situation we can define a map $h_{S}^{R} \colon
\pi_{R}(X) \to X_{S} \times \ker
\left( \pi_{S}^{R}|\pi_{R}(X)\right)$ by letting

\[ h_{S}^{R}(x) = \left( \pi_{S}^{R}(x),
x-r_{R}\left(\pi_{S}^{R}(x)\right) \right) \;\;
\text{for each}\;\; x \in \pi_{R}(X) .\]

A straightforward verification shows that $h_{S}^{R}$ is a
continuous linear map which satisfies the required equality
$\pi_{1}\circ h_{S}^{R} = \pi_{S}^{R}|\pi_{R}(X)$. Also note that by letting

\[ g_{S}^{R}(y,x) = r_{R}(y,0) + x \;\; \text{for each}\;\;
(y,x) \in \pi_{S}(X) \times \ker\left( \pi_{S}^{R}|\pi_{R}(X)\right) \]

\noindent we define a continuous linear map $g_{S}^{R} \colon
\pi_{S}(X) \times \ker\left( \pi_{S}^{R}|\pi_{R}(X)\right) \to \pi_{R}(X)$.
It is easy to see that 

\[g_{S}^{R}\circ h_{S}^{R} = \operatorname{id}_{\pi_{R}(X)}\;\;\text{and}\;\;
h_{S}^{R}\circ g_{S}^{R} = \operatorname{id}_{\pi_{S}(X) \times \ker\left( \pi_{S}^{R}|\pi_{R}(X)\right)} .\]

This proves that $h_{S}^{R}$ is a topological isomorphism 
and finishes the proof of Claim 3.
\medskip

{\bf Claim 4.} {\em Every countable subset of $T$ is
contained in a countable $r$-admissible subset of $T$.}
\medskip

Let $A$ be a countable subset of $T$. Our goal is to find a countable
$r$-admissible subset $S$ such that $A \subseteq S$. By Proposition
\ref{P:spectral}, there exist a countable subset $S$ of $T$
and a continuous homomorphism $r_{S} \colon B_{S} \to B_{S}$ such
that $A \subseteq S$ and $\pi_{S}\circ r = r_{S}\circ \pi_{S}$.
Consider a point $y \in \pi_{S}(X)$. Also pick a point $x \in X$
such that $\pi_{S}(x) = y$. Then
\[ y  = \pi_{S}(x) = \pi_{S}\left( r(x)\right) =
r_{S}\left(\pi_{S}(x)\right) = r_{S}(y) .\]
This shows that $r_{S}|\pi_{S}(X) = \operatorname{id}_{\pi_{S}(X)}$ 
(this shows, in fact, that $\pi_{S}(X)$ is complemented in $B_{S}$).

In order to show that $S$ is $r$-admissible let us consider a point
$z \in \pi_{S}^{-1}\left(\pi_{S}(X)\right)$. By the observation
made above,
$r_{S}\left(\pi_{S}(z)\right) = \pi_{S}(z)$. Finally
\[ \pi_{S}(z) = r_{S}\left(\pi_{S}(z)\right) =
\pi_{S}\left( r(x)\right) \]
which implies that $S$ is $r$-admissible.

We now use the above listed properties of $r$-admissible
subsets and proceed as follows.
By Claim 4, each element
$t_{\alpha} \in T$ is contained in a countable $r$-admissible
subset $S_{\alpha} \subseteq T$. According to Claim 1, the set
$T_{\alpha} = \bigcup\{ S_{\beta} \colon
\beta \leq \alpha\}$ is $r$-admissible for each $\alpha < \tau$.
Consider the projective system
\[ {\mathcal S}_{X} = \{ X_{\alpha}, p_{\alpha}^{\alpha +1}, \tau\} ,\]
where 
\[ X_{\alpha} = \pi_{T_{\alpha}}(X)\;\;\text{and}\;\;
p_{\alpha}^{\alpha +1} = \pi_{T_{\alpha}}^{T_{\alpha +1}}|
\pi_{T_{\alpha +1}}(X) \colon X_{\alpha +1} \to X_{\alpha}\;\;
\text{for each}\;\; \alpha < \tau .\]
Since $T = \bigcup\{ T_{\alpha} \colon \alpha < \tau\}$, it follows
that $X = \projlim\mathcal S$. Obvious transfinite induction based on
Claim 3 shows that 
\[ X = \projlim{\mathcal S} = X_{0} \times\prod\{
\operatorname{ker}\left( p_{\alpha}^{\alpha +1}\right) \colon
\alpha < \tau \} .\]
Since, by the construction, $S_{\alpha}$ is a countable
$r$-admissible subset of $T$, it follows from Claim 2 that
$X_{0}$ and $\ker\left( p_{\alpha}^{\alpha +1}\right)$, $\alpha <\tau$,
being complemented subspaces of countable products of Banach spaces, are
Fr\'{e}chet spaces. This finishes the proof of
Theorem \ref{T:complementedpro}.
\end{proof}

Recall that an object $X$ of the category
${\mathcal LSC}$ of locally convex spaces and their
continuous linear maps is {\em injective} if any continuous linear map
$f \colon A \to X$, defined on a linear subspace of a space $B$, admits
a continuous linear extension $g \colon B \to X$ (i.e. $g|A = f$).

The following statement is related to Problem 1 stated in the Introduction.

\begin{cor}\label{C:inj}
The following conditions are equivalent for a
locally convex topological vector space $X$:
\begin{itemize}
\item[(1)]
$X$ is an injective object of the category $\mathcal LCS$.
\item[(2)]
$X$ is isomorphic to the product
$\displaystyle \prod\{ F_{t} \colon t \in T\}$,
where each $F_{t}$, $t \in T$,
is a complemented subspace of a product
$\displaystyle \prod\{ \ell_{\infty}(J_{t_{n}}) \colon n \in \omega\}$.
\end{itemize}
\end{cor}
\begin{proof}
$(2) \Longrightarrow (1)$. By \cite[Lemma 0]{domort89} and \cite[p.105]{lintza77},
$\ell_{\infty}(J)$ is an injective
object of the category $\mathcal LCS$ for any set $J$. Obviously
(see, for instance, \cite[Lemma 1.9]{domort89}) product of an
arbitrary collection of injective objects of the category
$\mathcal LCS$ is also an injective
object of this category. Consequently, the Fr\'{e}chet space
$F_{t}$, $t \in T$, as a complemented subspace of
$\displaystyle \prod\{ \ell_{\infty}(J_{t_{n}})\colon n \in \omega\}$,
is injective. Finally, the space $X$, as a product of injectives,
is an injective object of the category $\mathcal LCS$.

$(1) \Longrightarrow (2)$. The space $X$ can be identified
with a closed linear subspace of the product
$\displaystyle \prod\{ B_{t} \colon t \in T\}$ of
Banach spaces $B_{t}$, $T \in T$. Each of the spaces $B_{t}$ can in turn
be identified with a closed linear subspace of the space
$\ell_{\infty}(J_{t})$ for some set $J_{t}$, $t \in T$.
Condition (1) implies in this situation that $X$
is a complemented subspace of the product 
$\displaystyle\prod\{ \ell_{\infty}(J_{t}) \colon t \in T\}$.
The required conclusion now follows from Theorem
\ref{T:complementedpro}.
\end{proof}


\providecommand{\bysame}{\leavevmode\hbox to3em{\hrulefill}\thinspace}



\end{document}